\documentclass[11pt]{amsart}
\usepackage{amsmath,amssymb,amscd}
\usepackage{graphicx}
\usepackage{a4wide}

\newcommand{\ts}{\hspace{0.5pt}}
\newcommand{\nts}{\hspace{-0.5pt}}
\newcommand{\oplam}{\mbox{\Large $\curlywedge$}}

\DeclareMathOperator{\cov}{cvg}
\DeclareMathOperator{\dens}{dens}
\DeclareMathOperator{\vol}{vol}

\begin{document}

\title{Homometric Point Sets and Inverse Problems}

\author{Uwe Grimm}
\address{Department of Mathematics and Statistics, The Open University,
Walton Hall,\newline\hspace*{\parindent}Milton Keynes MK7 6AA, 
United Kingdom}
\email{u.g.grimm@open.ac.uk}
\urladdr{http://mcs.open.ac.uk/ugg2/}

\author{Michael Baake}
\address{Fakult\"at f\"ur Mathematik, Universit\"at Bielefeld,
Postfach 100131, 33501 Bielefeld,  Germany}
\email{\texttt{mbaake@math.uni-bielefeld.de}}
\urladdr{\texttt{http://www.math.uni-bielefeld.de/baake}}

\begin{abstract}
The inverse problem of diffraction theory in essence amounts
to the reconstruction of the atomic positions of a solid from its
diffraction image. From a mathematical perspective, this is a
notoriously difficult problem, even in the idealised situation of
perfect diffraction from an infinite structure.

Here, the problem is analysed via the autocorrelation measure of the
underlying point set, where two point sets are called homometric when
they share the same autocorrelation. For the class of mathematical
quasicrystals within a given cut and project scheme, the homometry
problem becomes equivalent to Matheron's covariogram problem, in the
sense of determining the window from its covariogram. Although certain
uniqueness results are known for convex windows, interesting examples
of distinct homometric model sets already emerge in the plane.

The uncertainty level increases in the presence of diffuse
scattering. Already in one dimension, a mixed spectrum can be
compatible with structures of different entropy. We expand on this
example by constructing a family of mixed systems with fixed
diffraction image but varying entropy. We also outline how this
generalises to higher dimension.
\end{abstract}

\maketitle

\section{Introduction}

After 25 years of quasicrystal research, our understanding of the
atomic structure of quasicrystalline alloys is still far from being
complete \cite{Steurer}. The main reason for this is the difficult
inverse problem of determining the structure at the atomic scale from
the available information, which exists mainly in the form of
diffraction intensities.  Here, we discuss the non-uniqueness arising
from \emph{homometric} point sets, which is even present in the
idealised situation of a perfect diffraction measurement from an
infinite point set $\varLambda\subset\mathbb{R}^{d}$.

First, we consider the situation where $\varLambda$ is a mathematical
quasicrystal or model set. A perfect diffraction image of
$\varLambda$, as described by the positive diffraction measure
$\!\widehat{\,\gamma^{}_{\!\varLambda}}$, uniquely determines its inverse
Fourier transform, which is the autocorrelation (or Patterson) measure
$\gamma^{}_{\!\varLambda}$.  The starting point is thus the
(hypothetically complete) knowledge of $\gamma^{}_{\!\varLambda}$, and
the remaining task is then to determine $\varLambda$ from this
information.  For a model set based on a known cut and project scheme,
this amounts to determine the corresponding window $W$ in internal
space.

Beyond pure point diffraction, we reconsider the known homometry
between the binary Rudin-Shapiro sequence and the Bernoulli (or coin
flipping) chain \cite{HB00}. We introduce a new process, called
`Bernoullisation', which provides a continuous isospectral transition
between these two extremal cases. This method generalises to arbitrary
dimension and shows that even a perfect diffraction image (of mixed
type) may not be able to distinguish structures of different entropy.

\section{Homometry}

For finite point sets $F\subset\mathbb{R}^{d}$, homometry is defined
in terms of their \emph{difference sets} $F-F$, taking into account
multiplicities.  Two finite point sets are called \emph{homometric}
when they share the same weighted difference set (which is a
multi-set), meaning that each difference vector occurs with the same
cardinality in either set; see \cite{Patt44} for an early class of
examples in one dimension.

A relatively simple homometric pair, realised as finite subsets
$F_{1}\ne F_{2}\subset\mathbb{Z}^{2}$, was constructed in
\cite{GGZ05}. One choice of the coordinates results in
\begin{equation}\label{eq:f1f2}
\begin{split}
F_{1} =& \,\textstyle\bigl\{
        \binom{0}{0}, \binom{1}{0}, \binom{1}{1}, \binom{1}{2}, 
        \binom{1}{3}, \binom{2}{1}, \binom{2}{2},
        \binom{2}{3}, \binom{2}{4}, 
        \binom{2}{5}, \binom{3}{3}, \binom{3}{4}, 
        \binom{3}{5}, \binom{4}{4}, \binom{4}{5}
        \bigr\} ,\\
F_{2} =& \,\textstyle\bigl\{
        \binom{0}{0}, \binom{0}{1}, \binom{1}{0}, \binom{1}{1}, 
        \binom{1}{2}, \binom{1}{3}, \binom{1}{4}, 
        \binom{2}{2}, \binom{2}{3},
        \binom{2}{4}, \binom{2}{5}, \binom{3}{3}, 
        \binom{3}{4}, \binom{3}{5}, \binom{4}{5}
        \bigr\} .
\end{split}
\end{equation}
One can explicitly check that $F_{1}-F_{1}=F_{2}-F_{2}$, including
multiplicities.

An appropriate generalisation to infinite point sets needs the concept
of density.  We call two infinite point sets \emph{homometric} when
their natural autocorrelation measures exist and coincide.
Homometric point sets thus have the same density. Due to the
volume averaging involved, two point sets related by adding or
removing a point set of density $0$ are homometric. This is also true
of point sets related by translation or inversion (but not, in
general, by rotation).

It is well known that two crystallographic (or fully periodic) point
sets can only be homometric when they share the same lattice of
periods. They are then mutually locally derivable (MLD) from each
other \cite{BSJ,B}, which also implies that the associated dynamical
systems (under the translation action \cite{M}) are topologically
conjugate \cite{K}. The corresponding question for mathematical
quasicrystals (model sets without any periodicity) is more difficult,
as we shall demonstrate by an example.

\subsection{Covariogram}

There is an interesting connection between the homometry of model sets
(with a Euclidean internal space) and the covariogram problem. For a
non-empty, relatively compact subset $K\subset\mathbb{R}^{d}$, which is
assumed to be Riemann measurable, the function
\begin{equation}\label{eq:cov}
   \cov^{}_{K} (x)  := \vol\bigl(K\cap(x+K)\bigr),
\end{equation}
defined for all $x\in\mathbb{R}^{d}$, is called the \emph{covariogram}
of $K$.  The covariogram problem amounts to determine $K$ from its
covariogram $\cov^{}_{K}(x)$; compare \cite{Bia05,GGZ05}. This is
sometimes also referred to as Matheron's problem, which was originally
formulated as the question whether the covariogram determines a convex
body, among all convex bodies, up to translation and inversion; see
\cite{Mat75,Mat86,Bia05} for details. Since $\cov^{}_ {K}
(x)=\cov^{}_{t+K} (x)$ for any translation $t\in\mathbb{R}^{d}$ and
$\cov^{}_ {K} (x)=\cov^{}_{-K} (x)$, the covariogram $\cov^{}_{K}$ can
determine $K$ at best up to translations and inversion. We call two
non-empty, relatively compact, Riemann measurable sets
$K,K^{\prime}\subset\mathbb{R}^{d}$ \emph{homometric} when
$\cov^{}_{K}=\cov^{}_{K^{\prime}}$.

Denoting the characteristic function of $K$ by $1^{}_{K}$, the
function $\cov^{}_{K}$ is given by the convolution
\begin{equation} \label{eq:covconv}
   \cov^{}_{K} (x)  =  
   \bigl(1^{}_{K} * 1^{}_{-K}\bigr) (x) \ts .
\end{equation}
Its Fourier transform 
\begin{equation} \label{eq:covft}
   \widehat{\,\cov^{}_{K}} (k)  = 
   \big\lvert \widehat{\, 1^{}_{K}} (k) \big\rvert^2 
\end{equation}
is an analytic, positive function that vanishes in the limit as 
$\lvert k \rvert \to \infty$.  This
relation is the reason why, if $K$ is itself inversion symmetric in
the sense that $-K=t+K$ for a suitable translation
$t\in\mathbb{R}^{d}$, the function $1^{}_{K}$ (and hence $K$) can be
reconstructed from the knowledge of $\cov^{}_{K}$, up to translation
and inversion \cite{CJ94}.

If $K$ is a convex polytope in dimension $d\le 3$, it is determined by
$\cov^{}_{K}$; see \cite{Bia02,Bia05,AB07,Bia08} and references
therein. In general, however, the reconstruction of $K$ from the
knowledge of $\cov^{}_{K}$ is a difficult problem. An interesting
example of two polyominoes with the same covariogram \cite{BG07}
follows from the point set pair of Eq.~\eqref{eq:f1f2} by adding the
unit square $C=\bigl[-\frac{1}{2},\frac{1}{2}\bigr]^2$, so that
\begin{equation}\label{eq:p12}
   P_{1}=F_{1}+C \quad \text{and}\quad  P_{2}=F_{2}+C . 
\end{equation}
Their covariograms are equal as a consequence of the homometry of the
finite point sets $F_{1}$ and $F_{2}$, whence $P_{1}$ and $P_{2}$ are
homometric (as are also any translates of $\pm P_{1}$ and $\pm
P_{2}$). The polyominoes $P_{1}$, $P_{2}$ and their joint covariogram
are displayed and discussed in more detail in \cite{BG07}.  Let us
mention in addition that the scaled polyominoes $\alpha P_{1}$ and $\alpha
P_{2}$ are homometric to each other for any choice of
$\alpha\in\mathbb{R}$.

\subsection{Homometry of model sets}

Let us now consider the situation of regular model sets $\varLambda$
that are defined via a cut and project scheme \cite{Moody00,B} with
Euclidean internal space
\begin{equation}\label{eq:cps}
\renewcommand{\arraystretch}{1.2}\begin{array}{r@{}ccccc@{}l}
   & \mathbb{R}^{d} & \xleftarrow{\,\;\;\pi\;\;\,} & 
   \mathbb{R}^{d} \times \mathbb{R}^{m} & 
        \xrightarrow{\;\pi^{}_{\mathrm{int}\;}} & \mathbb{R}^{m} & \\
   & \cup & & \cup & & \cup & \hspace*{-1ex} 
   \raisebox{1pt}{\text{\footnotesize dense}} \\
   & \pi(\mathcal{L}) & \xleftarrow{\; 1-1 \;} & \mathcal{L} & 
   \xrightarrow{\; \hphantom{1-1} \;} & \pi^{}_{\mathrm{int}}(\mathcal{L}) & \\
   & \| & & & & \| & \\
   & L & \multicolumn{3}{c}{\xrightarrow{\qquad\qquad\quad\star
    \qquad\qquad\quad}} 
       &  {L_{}}^{\star\nts} & \\
\end{array}\renewcommand{\arraystretch}{1}
\end{equation}
by
\begin{equation}\label{eq:modelset}
    \varLambda=t+\oplam(W)=
    t+\bigl\{x\in\pi(\mathcal{L})\mid x^{\star}\in W\bigr\} .
\end{equation}
Here, $\mathcal{L}$ is a lattice in
$\mathbb{R}^{d}\times\mathbb{R}^{m}$, the window
$W\subset\mathbb{R}^{m}$ is a non-empty, relatively compact set with
boundary of measure $0$, and $t\in\mathbb{R}^{d}$ is an
arbitrary translation.
 
The autocorrelation $\gamma^{}_{\!\varLambda}$ of the corresponding
Dirac comb $\delta_\varLambda=\sum_{x\in\varLambda}\delta_{x}$ exists
and has the explicit form $\gamma^{}_{\!\varLambda} =
\sum_{z\in\varLambda-\varLambda} \eta(z)\, \delta_z$, with
coefficients
\begin{equation} \label{eq:autocoeffwin}
   \eta(z) \,  =  \; \dens(\varLambda)\,
   \frac{\vol \bigl(W\cap(W-z^{\star\nts})\bigr)}{\vol (W)} \,
    =  \; \dens(\mathcal{L})\, \cov^{}_{W}(z^{\star\nts}) \ts ,
\end{equation}
expressed in terms of the covariogram of the window $W$. Hence, two
(Euclidean) model sets obtained from the same cut and project scheme
are homometric if and only if the defining windows share the same
covariogram.  As a consequence, homometric model sets from the same
cut and project scheme have the same diffraction measure. Conversely,
kinematic diffraction cannot discriminate between homometric model
sets.

A planar example is obtained by using the homometric pair of
polyominoes $P_{1}$ and $P_{2}$ as windows for model sets in a cut and
project scheme of type \eqref{eq:cps}. For instance, we can use the
Minkowski embedding
$\mathcal{L}_{8}\subset\mathbb{C}^{2}\simeq\mathbb{R}^{4}$ of
$L=\mathbb{Z}[\xi_{8}]$, where $\xi_{8}$ is a primitive $8$th root of
unity, and a $\star$-map defined by a suitable algebraic conjugation.
The two model sets $\varLambda_{1} := \oplam(P_{1})$ and
$\varLambda_{2} := \oplam(P_{2})$, with $P_{1}$ and $P_{2}$ as defined
in \eqref{eq:p12}, are then homometric by construction (the relative
position of the windows, which is irrelevant for homometry, maximises
their intersection).  The two model sets $\varLambda_{1}$ and
$\varLambda_{2}$ are \emph{not} locally indistinguishable, and differ
in points of positive density. In particular, the difference sets
$\varLambda_{1}\setminus\varLambda_{2}$ and
$\varLambda_{2}\setminus\varLambda_{1}$ are model sets themselves (but
not homometric).

The diffraction measure $\widehat{\gamma}$ is the same for both
$\varLambda_{1}$ and $\varLambda_{2}$, and reads
\[
   \widehat{\gamma}\: = 
   \sum_{k\in\frac{1}{2}L}  I(k^{\star\nts})\, \delta_{k} \ts ,
\]
with intensity function $I(y)=\lvert A_{i}(y)\rvert^{2}$ derived from
\[
   A_{i}(y)  = 
    \,\dens (\mathcal{L}^{}_{8})\, \widehat{\, 1^{}_{P_i}}(-y) \ts .
\]
While the amplitudes depend on the window, their absolute squares do
not.  One can work out the explicit diffraction intensities; see
\cite{BG07} for details. The ratio of the (complex) amplitudes is
given by
\begin{equation}
   \frac{A^{}_{1}(y)}{A^{}_{2}(y)} \,=\, 
   \frac{1+e^{2\pi i y^{}_{2}} + e^{2\pi i (y^{}_{1}+2y^{}_{2})}}{1 + 
         2\,  e^{\pi i (2y^{}_{1}+3y^{}_{2})}\cos(\pi y^{}_{2})}
   \,=\, 1  +  \frac{1-e^{2\pi i y^{}_{1}}}{e^{2\pi i y^{}_{1}} + 
       e^{2\pi i (y^{}_{1}+y^{}_{2})} + e^{-2\pi i y^{}_{2}}}
\end{equation}
with $y=(y^{}_{1},y^{}_{2})$. This is a well-defined function on
internal space $\mathbb{R}^{2}$, with values in $\mathbb{S}^{1}$,
unless the denominator vanishes. The latter happens for
$y^{}_{2}\in\mathbb{Z}+\{\frac{1}{3},\frac{2}{3}\}$ together with
$y^{}_{1}\in\mathbb{Z}$. One can check that the ratio has no
continuous extension to these points. Writing the ratio as
$\exp\bigl(2\pi i\ts\chi(y)\bigr)$, the phase function $\chi$ is not
defined at these points. Moreover, as one can check explicitly, it
does \emph{not} satisfy the additivity property
$\chi(y+y')=\chi(y)+\chi(y') \bmod 1$, wherefore the ratio fails to be
a character on $\mathbb{R}^{2}$ by violating both defining properties;
compare \cite{M} and references therein. Note that an analogous
phenomenon already shows up in the comparison of $P_{i}$ with
$-P_{i}$, because these windows are not inversion symmetric up to
translations.

The choice of the windows $P_{1}$ and $P_{2}$ is special in the sense
that $\varLambda_{1}$ and $\varLambda_{2}$ turn out to be MLD, because
the square $C$ satisfies 
\[
   C = P_{1} \cap (-t+P_{1}) = P_{2} \cap (-t+P_{2})
\]
with the translation $t=(4,5)$, and each window is now the union of
$15$ integral (and hence admissible) translates of $C$ according to
Eq.~\eqref{eq:p12}; see \cite{BSJ} for details. In this case, the
associated dynamical systems \cite{M} are again topologically
conjugate. As mentioned above, the two windows may be scaled (by the
same factor) without affecting their mutual homometry. For almost all
choices of the scaling factor, one loses the MLD property of the
corresponding model sets, because the finite reconstruction property
\cite{BSJ} is lost. Nevertheless, the associated dynamical systems
will always be \emph{metrically} isomorphic (due to the Halmos-von
Neumann theorem).  It is an interesting open question whether they are
still also \emph{topologically} conjugate, which is a weaker
equivalence notion than MLD.

Independent of this conjugacy issues, our example illustrates that
diffraction (hence autocorrelation) alone is generally insufficient to
uniquely determine a regular model set. However, as discussed in
\cite{M}, this ambiguity can be resolved with the knowledge of the
$3$-point correlations. This statement is immediate in our example
(via the existence or non-existence of certain patches), but holds in
full generality for regular model sets; see \cite{M} and references
therein.

\section{Random Dirac combs}

The problem of reconstruction becomes even more involved in the case
of mixed spectra. In this setting, a slight change in point of view is
helpful to separate distinct spectral components. This is most easily
achieved by considering \emph{weighted} Dirac combs of point sets,
with real (or even complex) weights. Below, generalising an example
discussed in \cite{HB00}, we construct a family of one-dimensional
homometric (weighted) point sets, based on the binary Rudin-Shapiro
sequence, which cover the entire entropy range from $0$ to $\log(2)$,
the maximal possible entropy for a binary system. This shows that, in
general, it is not even possible to determine the degree of long-range
order of the weighted point set from diffraction data. The same
conclusion also holds for the diffraction of the associated unweighted
point sets.

\subsection{Bernoulli versus Rudin-Shapiro}

We start by re-considering the example of Ref.~\cite{HB00}. The first
model is a Bernoulli system on $\mathbb{Z}$, with the stochastic Dirac
comb
\begin{equation}\label{eq:bern}
   \omega^{}_{\mathrm{B}}=\sum_{m\in\mathbb{Z}} 
    Y^{}_{\! m}\ts\delta^{}_{m} \ts ,
\end{equation}
where $(Y^{}_{\! m})^{}_{m\in\mathbb{Z}}$ is a family of
i.i.d.\ random variables that each take the values $1$ and $-1$, with
probabilities $p$ and $1-p$, where $0\le p \le 1$. 
For the stochastic Dirac comb $\omega^{}_{\mathrm{B}}$,
the autocorrelation  measure  $\gamma^{}_{\mathrm{B}}$
and the diffraction measure $\!\widehat{\,\gamma^{}_{\mathrm{B}}\!}$
almost surely exist and read
\begin{equation}\label{eq:diffbern}
   \begin{split}
   \gamma^{}_{\mathrm{B}} &= \, (2p-1)^{2}\delta^{}_{\mathbb{Z}} + 
   4\ts p\ts (1-p)\,\delta^{}_{0}\ts ,\\
   \widehat{\,\gamma^{}_{\mathrm{B}}\!}\, &= 
    \, (2p-1)^{2}\delta^{}_{\mathbb{Z}} + 
   4\ts p\ts (1-p)\,\lambda\ts ,
   \end{split}
\end{equation}
where $\lambda$ denotes Lebesgue measure on $\mathbb{R}$. Note that,
in this stochastic situation, almost sure results are unavoidable.  In
particular, one has $\!\widehat{\,\gamma^{}_{\mathrm{B}}\!}\,=\lambda$
for $p=\frac{1}{2}$ and
$\!\widehat{\,\gamma^{}_{\mathrm{B}}\!}\,=\delta^{}_{\mathbb{Z}}$ for
$p=0$ or $p=1$. The choices $p=0$ and $p=1$ correspond to the
deterministic limiting cases
$\omega^{}_{\mathrm{B}}=\pm\ts\delta^{}_{\mathbb{Z}}$, while
$p=\frac{1}{2}$ describes a stochastic comb (coin tossing) with
weights of average $0$.

The binary Rudin-Shapiro sequence is defined in two steps as follows
\cite{Q}. We start from the substitution rule
\begin{equation}\label{eq:rs}
    a\mapsto ac\ts, \quad b\mapsto dc\ts,\quad
    c\mapsto ab\ts,\quad d\mapsto da\ts,
\end{equation}
on the four-letter alphabet $\mathcal{A}=\{a,b,c,d\}$. We choose a
bi-infinite fixed point (under the square of the above substitution,
with seed $ba$) and apply the morphism $\varphi\!: \mathcal{A}
\longrightarrow \{\pm 1\}$ defined by $\varphi(a)=\varphi(c)=1$ and
$\varphi(b)=\varphi(d)=-1$, extended to $\mathcal{A}^{\mathbb{Z}}$.
The autocorrelation and diffraction measures of the resulting binary
Rudin-Shapiro chain $S_{\mathrm{RS}}$ are
\begin{equation}\label{eq:rsdiff}
   \gamma^{}_{\mathrm{RS}} = \delta^{}_{0}
   \quad\text{and}\quad 
   \widehat{\gamma^{}_{\mathrm{RS}}} = \lambda\ts .
\end{equation}
This is an example with a purely absolutely continuous diffraction,
despite the fact that the Rudin-Shapiro chain is deterministic and has
entropy $0$. In particular, it agrees with the diffraction measure of
the Bernoulli comb with $p=\frac{1}{2}$, which has entropy $\log(2)$.

\subsection{`Bernoullisation'}

It is possible to impose the influence of chance on the order of a
deterministic system, and thus interpolate between deterministic and
random systems.  Here, we focus on binary sequences and modify them by
an i.i.d.\ family of Bernoulli variables.

Consider a bi-infinite binary sequence $S\in\{\pm 1\}^{\mathbb{Z}}$ which we
assume to be uniquely ergodic (in the sense that its hull under the
action of the shift map is a uniquely ergodic dynamical system). Then, the
corresponding Dirac comb $\omega^{}_{S} = \sum_{i\in\mathbb{Z}}
S_{i}\,\delta_{i}$ possesses the (natural) autocorrelation
$\gamma^{}_{S}=\sum_{m\in\mathbb{Z}}\eta^{}_{S}(m)\,\delta_{m}$ with 
autocorrelation coefficients $\eta^{}_{S}(m)$, where $\eta^{}_{S}(0)=1$ by
construction. 

Let $(Y^{}_{i})^{}_{i\in\mathbb{Z}}$ be an i.i.d.\ family of random variables
that each take values $+1$ and $-1$ with probabilities $p$ and $1-p$.
The \emph{`Bernoullisation'} of $\omega^{}_{S}$ is the random Dirac comb
\begin{equation}\label{eq:bernoullisation}
     \omega^{}_{S;p}\, := \sum_{i\in\mathbb{Z}} S_{i}\,Y_{\! i}\,\delta_{i} \ts ,
\end{equation}
which emerges from $\omega^{}_{S}$ by independently changing the sign
of each $S_{i}$ with probability $1-p$. Setting $Z_{i}:=S_{i}Y_{\! i}$
defines a new family of independent (though, in general, not
identically distributed) random variables, with values in $\{\pm
1\}$. Despite this modification, the autocorrelation $\gamma^{}_{S;p}$
of $\omega^{}_{S;p}$ almost surely exists and can be determined via
its autocorrelation coefficients $\eta^{}_{S;p}(m)$ as follows. Since
one always has $\eta^{}_{S;p}(0)=\eta^{}_{S}(0)=1$, let $m\ne 0$ and
consider, for large $N$, the sum
\begin{equation}\label{eq:coeff}
    \frac{1}{2N\!+\!1} \sum_{i=-N}^{N} Z_{i}\ts Z_{i-m}
    \; = \, \frac{1}{2N\!+\!1} \Bigl( \sum_{(+,+)} + \sum_{(-,-)} - 
        \sum_{(+,-)} - \sum_{(-,+)}
       \Bigr) Y_{\! i}\, Y_{\! i-m}\ts ,
\end{equation}
which is split according to the value of $(S_{i},S_{i-m})$. Each of
the four sums can be handled in the same way as for the Bernoulli
comb, thus contributing $(2p-1)^2$ times the frequency of the
corresponding sign pair.  Observing that the overall signs are the
products $S_{i}\ts S_{i-m}$, it is clear that, as $N\to\infty$, one
obtains (almost surely)
\begin{equation}
    \eta^{}_{S;p}(m) = (2p-1)^{2}\ts\eta^{}_{S}(m)
\end{equation}
for all $m\ne 0$. Thus, the autocorrelation $\gamma^{}_{S;p}$ of
$\omega^{}_{S;p}$ almost surely exists and is given by
\begin{equation}
   \gamma^{}_{S;p} = (2p-1)^{2}\ts\gamma^{}_{S} + 
   4\ts p\ts (1-p)\, \delta^{}_{0}\ts ,
\end{equation}
where $\gamma^{}_{S}$ is the unique autocorrelation of
$\omega^{}_{S}$.

Consider now the Bernoullisation (with parameter $p$) of the binary
Rudin-Shapiro sequence with random Dirac comb
$\omega^{}_{\mathrm{RS};p}$. Its autocorrelation measure almost surely
exists and reads $\gamma^{}_{\mathrm{RS};p}=\delta^{}_{0}$,
independently of $p$. This means that the random Dirac combs
$\omega^{}_{\mathrm{RS};p}$, even for different values of $p$, are
almost surely homometric, and share the purely absolutely continuous
diffraction measure $\widehat{\gamma^{}_{\mathrm{RS};p}}=\lambda$.

\subsection{Extension to higher dimension}

Our above discussion has an immediate extension to Euclidean space of
arbitrary dimension $d$. Consider $d$ complex-valued sequences
$(U^{(\ell)}_{i})^{}_{i\in\mathbb{Z}}$ and define the weighted Dirac comb
\[
    \omega^{}_{U}\, = 
     \, U^{(1)}\delta^{}_{\mathbb{Z}}\otimes\ldots\otimes
    U^{(d)}\delta^{}_{\mathbb{Z}} \, = 
    \sum_{x\in\mathbb{Z}^{d}} 
    \Bigl( \prod_{\ell=1}^{d} U^{(\ell)}_{x^{}_{\ell}} \Bigr)\,
    \delta_{x}\ts .
 \]
Here, $x=(x^{}_{1},\ldots,x^{}_{d})$, and the weights (on
$\mathbb{Z}^{d}$) are products of $d$ elements of the individual
sequences.  Assuming that the natural autocorrelations of the
individual sequences exist, the relevant observation is that the
resulting autocorrelation of $\omega^{}_{U}$ (and hence also the
corresponding diffraction measure) factorises accordingly.

Each $U^{(\ell)}$ may be chosen as a member of our previous
one-parameter family of \eqref{eq:bernoullisation}, in particular as
$U^{(\ell)}=S_{\mathrm{RS}}$ for all $\ell$. This results in a deterministic
weighted Dirac comb on $\mathbb{Z}^{d}$ with diffraction measure
$\widehat{\gamma}=\lambda$, where $\lambda$ now denotes Lebesgue measure
on $\mathbb{R}^{d}$. This represents a system of entropy $0$.

The analogue of the Bernoullisation of \eqref{eq:bernoullisation} in
$d$-space, with an i.i.d.\ family $(Y_{\! x})_{x\in\mathbb{Z}^{d}}$,
then results in an isospectral one-parameter family of random Dirac
combs on $\mathbb{Z}^{d}$, which can realise any entropy between $0$
and $\log(2)$. Similarly, using $\mathbb{Z}^{d}\simeq
\mathbb{Z}^{k}\times\mathbb{Z}^{d-k}$ for some $0\le k\le d$ and
restricting the Bernoullisation to $\mathbb{Z}^{k}$, one obtains
isospectral families with arbitrary entropy of rank $k$ between $0$
and $\log(2)$. This indicates that the variety of homometric
structures grows with dimension.

\section{Conclusions}

The homometry problem for regular model sets in dimensions $d\le 3$
appears to have a unique solution if one may assume that the window is
convex \cite{AB07,Bia08}. However, this favourable situation is not
always met in real quasicrystals. Our example, with non-convex
windows, illustrates the existence of distinct homometric structures.
Since homometric model sets are always metrically isomorphic, but not
necessarily MLD, it remains an interesting question how they are
related as topological dynamical systems.

For the case of spectra with a continuous component, the
Bernoullisation approach can explore the full entropy range: the
Bernoulli case (with $p=\frac{1}{2}$) has entropy $\log 2$, the
maximal value for a binary system, while Rudin-Shapiro has entropy
$0$, and the parameter $p$ interpolates continuously between the two
limiting cases. This gives an indication of how degenerate the inverse
problem really is.  Unless additional information is available, one
possible strategy to proceed would employ an optimisation approach,
for instance by choosing the structure which maximises the entropy,
which singles out the Bernoulli comb here. \bigskip

\section*{Acknowledgements}

We are grateful to R.V.~Moody for valuable discussions and comments on
the manuscript. It is a pleasure to thank the School of Mathematics
and Physics at the University of Tasmania for their kind hospitality
during an extended stay in Hobart.  This work was supported by the
German Research Council, within the Collaborative Research Centre 701,
and by EPSRC via Grant EP/D058465.

\end{document}